\documentclass{amsart} % article/book/letter/amsart
\usepackage[utf8]{inputenc}
\usepackage[T1]{fontenc}	% accents
\usepackage[french,english]{babel}
\usepackage{lmodern}			% police de caractères

\usepackage{graphicx}	% insertion des images / graphiques
\usepackage{placeins}
\usepackage{float} % pour mettre H dans les figures

\usepackage{color}
\definecolor{vertFonce}{rgb}{0,0.5,0}
\definecolor{numLignes}{rgb}{0.17,0.57,0.7}	%{43,145,175}
\definecolor{gris}{rgb}{0.5,0.5,0.5}
\definecolor{grisFonce}{rgb}{0.2,0.2,0.2}
\definecolor{orange}{rgb}{1,0.65,0.31}		%{255,167,79}
\definecolor{orangeFonce}{rgb}{1,0.4,0}
\definecolor{bleuFonce}{rgb}{0,0,0.4}
\definecolor{rougeFonce}{rgb}{0.3,0,0}
\definecolor{rougeWord}{rgb}{0.5,0,0}
\definecolor{vertClair}{rgb}{0.8,1,0.8}
\definecolor{rougeClair}{rgb}{1,0.5,0.5}

\usepackage{pict2e}
\usepackage{multido}

\usepackage{amsfonts,amssymb,amsthm,amsmath} 
\usepackage{dsfont}				% utilisation des caractères \mathds
\usepackage{mathrsfs}
\usepackage{yfonts}
 
\newtheorem{theorem}{Theorem}

\newtheorem{remark}{Remark}[section]

\newenvironment{demo}[1][]{%
	\begin{proof}[\textbf{Proof#1}]
	}{%
	\end{proof}
}

\newcommand{\step}[1]	{\paragraph{\itshape\bfseries Step #1.}}

\newcommand		{\R}		{\mathbb R}			% réels
\newcommand		{\DD}		{\mathcal{D}}

\newcommand			{\dd}		{\mathrm{d}}	% différentielle
\renewcommand		{\d}		{\,\dd}			% différentielle (intg)
\newcommand			{\dt}		{\frac{\dd}{\dd t}}
\DeclareMathOperator{\divg}		{div}			% divergence
\newcommand		{\indic}		{\mathds{1}}	% indicatrice

\newcommand		{\lt}			{\left}			%
\newcommand		{\rt}			{\right}		%
\newcommand		{\lal}			{\langle}		%
\newcommand		{\ral}			{\rangle}		%
\newcommand		{\weight}[1]	{\lal #1\ral}	% <x>
\newcommand		{\intd}			{\int_{\R^d}}
\newcommand		{\iintd}		{\iint_{\R^{2d}}}

\newcommand		{\init}			{\mathrm{in}}
\newcommand		{\loc}			{\mathrm{loc}}
\newcommand		{\eps}			{\varepsilon}

\newcommand		{\CS}			{\mathcal{C}^\mathrm{S}}
\newcommand		{\CHLS}			{\mathcal{C}^\mathrm{HLS}}

\newcommand		{\ak}			{a}
\newcommand		{\al}			{\alpha}

\usepackage{stmaryrd}
\newcommand		{\subsetArrow}	{\mathrel{\ooalign{$\subset$\cr%
\hidewidth\raise-.087ex\hbox{$_\shortrightarrow\mkern-1.5mu$}\cr}}}
\newcommand		{\subsetarrow}	{\mathrel{\ooalign{$\subset$\cr%
\hidewidth\raise-1.45ex\hbox{$\vec{}\mkern6mu$}\cr}}}

\title[p-Laplacian Keller-Segel Equation]{p-Laplacian Keller-Segel Equation: Fair Competition and Diffusion Dominated Cases}
\author{Laurent Lafleche and Samir Salem}

\date{\today}
 \def\signll{\bigskip\begin{center}{
 	\sc Laurent
 	Lafleche\par\vspace{2mm}
 	Université Paris-Dauphine, PSL Research University\par
 	CEREMADE, UMR CNRS 7534\par
 	Place du Maréchal de Lattre de Tassigny \par
 	75775 Paris Cedex 16 FRANCE\par
 	\vspace{2mm} e-mail:}
 	\tt{lafleche@ceremade.dauphine.fr}
 \end{center}}
 
 \def\signss{\bigskip\begin{center}{
 	\sc Samir
 	Salem\par\vspace{2mm} 
 	Université Paris-Dauphine, PSL Research University\par
 	CEREMADE, UMR CNRS 7534\par
 	Place du Maréchal de Lattre de Tassigny \par
 	75775 Paris Cedex 16 FRANCE\par
 	\vspace{2mm} e-mail:}
 	\tt{salem@ceremade.dauphine.fr}
 \end{center}}

\begin{document}

\date{\today}

\begin{abstract} 
This work deals with the aggregation diffusion equation	
	\begin{equation*}
		\partial_t \rho = \Delta_p\rho + \lambda \divg\lt((K_\ak*\rho)\rho\rt),
	\end{equation*}
	where $K_\ak(x)=\frac{x}{|x|^\ak}$ is an attraction kernel and $\Delta_p$ is the so called $p$-Laplacian. We show that the domain $\ak < p(d+1)-2d$ is subcritical with respect to the competition between the aggregation and diffusion by proving that there is existence unconditionally with respect to the mass. In the critical case we show existence of solution in a small mass regime for an $L\ln L$ initial condition.\newline
\textsc{Résumé.} Ce travail concerne l'étude d'une famille d'équations d'agrégation diffusion    
	\begin{equation*}
	\partial_t \rho = \Delta_p\rho + \lambda \divg\lt((K_\ak*\rho)\rho\rt),
	\end{equation*}
	où $K_\ak(x)=\frac{x}{|x|^\ak}$ est un champ d'attraction et $\Delta_p$ est le $p$-Laplacien. On montre que le domaine $\ak < p(d+1)-2d$  est sous-critique du point de vue de la compétiton entre l'agrégation et la diffusion en montrant l'existence de solution quelle que soit la masse. Dans le cas critique, on montre l'existence de solution dans un régime de petite masse pour une condition $L\ln L$.  

\end{abstract}

\maketitle

%%%%%%%%%%%%%%%%%%%%%%%%%%%%%%%%%%%%%%%%%%%%%%%%%%%%%%%%%%%%%%%%%%%%%%%%%%%%%%%%%
%
%
%                        Section: Introduction
%
%
%%%%%%%%%%%%%%%%%%%%%%%%%%%%%%%%%%%%%%%%%%%%%%%%%%%%%%%%%%%%%%%%%%%%%%%%%%%%%%%%%

\selectlanguage{french}
\section{Version francaise abrégée}

	On entend par équation d'agrégation-diffusion une equation aux dérivées partielles non linéaire sur $\R^d$ de la forme
	\begin{equation*}
		\partial_t \rho = \DD(\rho) + \lambda \divg\lt((K_\ak*\rho)\rho\rt),
	\end{equation*}
	où pour $\ak \in (0,d)$, $K_\ak(x)=\frac{x}{|x|^\ak}$ est un noyau d'attraction, $\lambda > 0$ indique l'intensité de cette interaction et $\DD$ est un opérateur de diffusion. Cette équation décrit par exemple l'évolution de la densité d'une population de bactéries ou d'astres en gravitation (voir par exemple \cite{hoffmann_keller-segel-type_2017}).
	
	Ce modèle a été largement étudié dans le cas de l'opérateur de diffusion nonlinéaire $\DD(\rho)=\Delta(\rho^m)$ pour $m>0$ (voir \cite{calvez_equilibria_2017}). Grâce à la structure algébrique conférée par ce choix de diffusion, on montre que l'EDP est en fait un flot de gradient pour la distance de Wasserstein d'ordre $2$ d'une certaine fonctionelle. De l'étude de cette fonctionnelle découle que la ligne $\ak = 2 - d(m-1)$ (dans le plan $(m,\ak)$) est critique du point de vue de la compétition entre l'agrégation et la diffusion. Le demi-plan situé au dessus de cette droite correspond au régime d'agrégation dominante, et celui au dessous à celui de diffusion dominante
	
	Lorsque la diffusion est fractionaire i.e. $\DD= \Delta^{\al/2}$ est le Laplacien fractionnaire d'exposant $\al\in(0,2)$, on montre que la ligne critique est la première bissectrice $\ak=\al$ (voir \cite{salem_propagation_2017,lafleche_fractional_2018}), et qu'elle délimite dans ce cas également deux régimes opposés. 
	
	Cette note poursuit cette étude, dans le cas du p-Laplacien $\DD(\rho) = \Delta_p(\rho) = \divg(|\nabla \rho|^{p-2}\nabla \rho)$ pour $p\in \lt(\frac{2d}{d+1},\frac{3d}{d+1}\rt)$. On montre ici que le domaine $\ak < p(d+1)-2d$ est sous-critique et qu'il y a existence pour petite masse dans le cas d'égalité. Au passage on établit une estimation de moments pour la $p$-équation de la chaleur.

\selectlanguage{english}
\section{Introduction}\label{intro}

	Aggregation diffusion equations play an important role in the modeling of collective behavior and more specially, in the case of the motion of cells and bacteria (see for instance \cite{hoffmann_keller-segel-type_2017}). The  (parabolic-elliptic) Keller-Segel equation, which has been extensively studied  (see \cite{blanchet_two-dimensional_2006}), is a typical example. In generality, we mean by aggregation equation the class of mean field nonlinear conservation equation of the form 
	\begin{equation}\label{eq:agdif}
		\partial_t \rho = \DD(\rho) + \lambda \divg\lt((K_\ak*\rho)\rho\rt),
	\end{equation}
	where $K_\ak$ is an aggregation kernel defined as $K_\ak(x) = \frac{x}{|x|^\ak}$, $\lambda>0$ is a parameter encoding the intensity of the aggregation and $\DD$ is some diffusion operator. Equation \eqref{eq:agdif} can then be interpreted as the evolution of the probability density of particles attracting each other through $K_\ak$ and diffusing through $\DD$. Then depending on the result of the competition between these two phenomena, the equation may yield to global existence or finite time blow up.
	
	The case of power law diffusion $\DD(\rho)=\Delta (\rho^m)$ for some $m>0$, has been studied in \cite{calvez_equilibria_2017} where the line $\ak= 2 - d(m-1)$ is shown to be critical. In that case equation \eqref{eq:agdif} can be seen as the gradient flow of some suitable functional with respect to the Wasserstein-$2$ distance and the criticality appears from the asymptotic study of this functional.
	
	The case of fractional diffusion $\DD(\rho) = \Delta^{\al/2}\rho$ for some $\al\in(0,2)$ has been studied in \cite{salem_propagation_2017,lafleche_fractional_2018}, where it is shown that the critical line is the first bisector $\al=\ak$, above which blow up of solutions may occur in finite time, and under which global well-posedness and propagation of chaos hold. 
	
	In order to complete this study, this note investigates the case where the diffusion operator is the $p$-Laplacian $\DD = \Delta_p$ (see e.g. \cite{lindqvist_notes_2006}), which is defined for any $\rho\in W^{1,p-1}_\loc$ by
	\begin{equation*}
		\forall\varphi\in C^\infty_c(\R^d),\ \lal \Delta_p \rho,\varphi \ral = -\intd \lt| \nabla \rho \rt|^{p-2}\nabla\rho \cdot\nabla \varphi,
	\end{equation*}
	and appears for example in the diffusion equations for sandpiles (see e.g. \cite{aronsson_fast/slow_1996,evans_fast/slow_1997}).

\section{Main results} 

	The aggregation equation \eqref{eq:agdif} with $\DD=\Delta_p$,
	\begin{equation}\label{eq:pLapKS}
		\partial_t\rho = \divg(|\nabla \rho|^{p-2}\nabla \rho) + \lambda \divg((K_\ak*\rho)\rho),
	\end{equation}
	has not been much studied, to the best of the author's knowledge. The only reference at this matter is \cite{liu_degenerate_2016}, which concerns the case $\ak=d$ and $p\in (2,\frac{3d}{d+1})$.
	
	\begin{figure}[H]\centering\label{fig1}
		\mbox{
			{\includegraphics[scale=0.7]{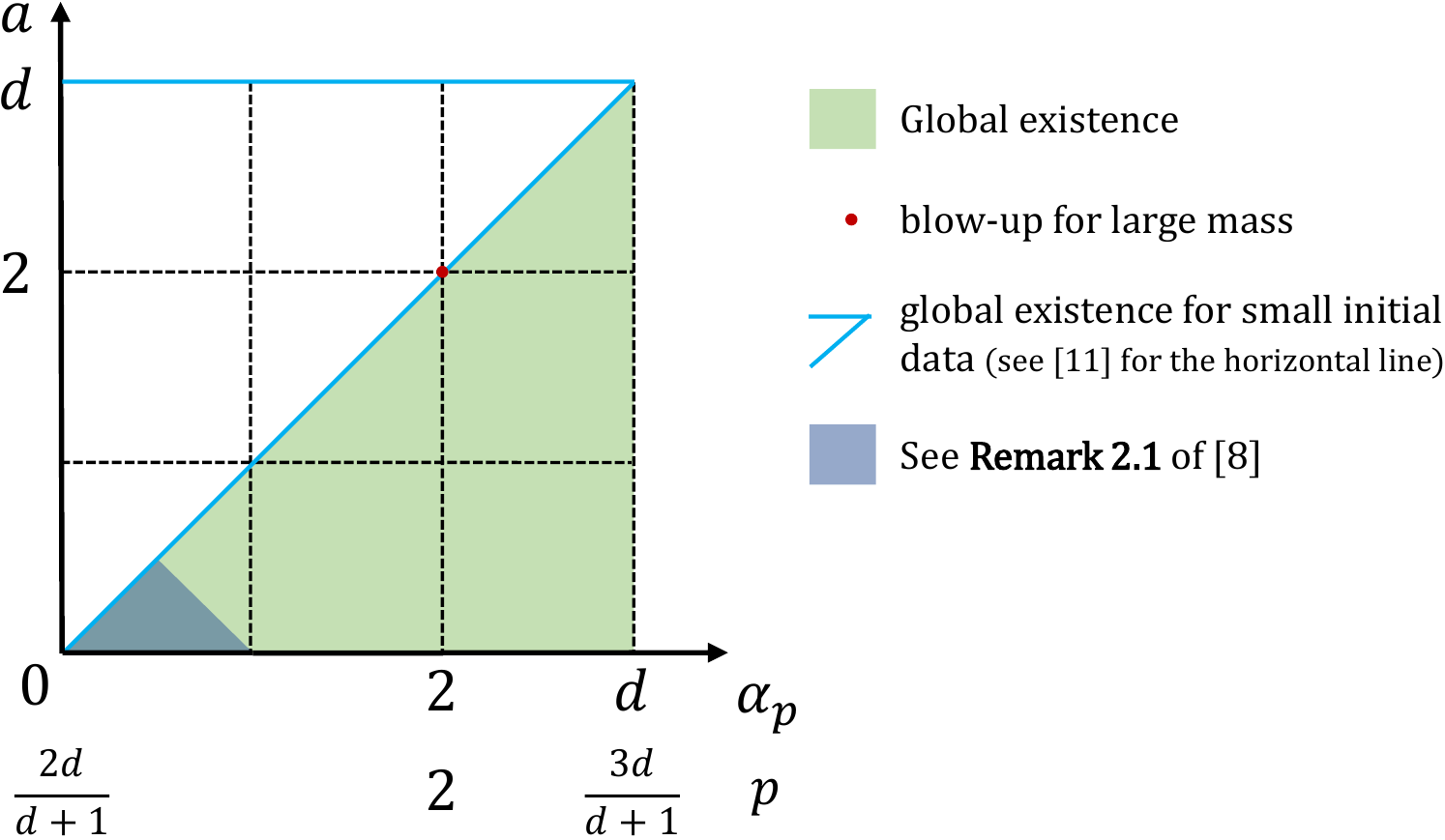}}}
		\caption{Graph of results about equation~\eqref{eq:pLapKS}} {{\small Graphique des résultats pour l'équation~\eqref{eq:pLapKS}}}
	\end{figure}
	
	Denoting $\|\rho\|_{L^p_k} := \|\rho m\|_{L^p}$ with $m(x)=\lal x\ral^k$ and $L\ln L = \{\rho \in L^1, \rho\ln \rho\in L^1\}$, we state the main result of this note.

	\begin{theorem}\label{thm:main}
		Let $d\geq2$, $\lambda>0$ and $(\ak,p)\in(0,d)\times\lt(\frac{2d}{d+1},\frac{3d}{d+1}\rt)$. Denote $\al_{p}:=p(d+1)-2d$ and assume $\al_p+\ak>1$. Let $\rho^\init\in L\ln L\cap L_k^1$ for some $k\in \lt( (1-\ak)_+,\alpha_p\wedge 1 \rt)$. Then in the
		\begin{itemize}
			\item \textbf{Diffusion dominated} case $\ak<\al_{p}$ 
			\item \textbf{Fair competition} case $\ak=\al_{p}$ if $\rho^\init$ satisfies
			\begin{align*}
				&M_0:=\|\rho^\init\|_{L^1}< C_{d,p} \lambda^{-\frac{1}{3-p}},
			\end{align*}
		\end{itemize}
		there exists a solution $\rho\in L^{p/p'}_\loc(\R_+,L^{p^*/p'})\cap L^\infty_\loc(\R_+,L^1_k)$ to equation \eqref{eq:pLapKS} with initial condition $\rho^\init$.
	\end{theorem}

	\begin{remark}
		The constant $C_{d,p}$ is given by
		\begin{equation*}
			C_{d,p} = \lt((d-\alpha_p) \CHLS_{d,\alpha_p,\frac{2d}{2d-\al_p}}\lt(p'^{-1}\CS_{d,p}\rt)^{p}\rt)^{-\frac{1}{3-p}},
		\end{equation*}
		where for $\ak\in(0,d)$ and $2-\frac{\ak}{d}=\frac{2}{q}$, $\CHLS_{d,\ak,q}$ is the best constant for the Hardy-Litllewood-Sobolev's inequality, 
		\begin{equation*}
			\iintd |x-y|^{-\ak}\rho(x)\rho(y)\d x\d y \leq C^{HLS}_{d,\ak,q}\|\rho\|^2_{L^q},
		\end{equation*}
		and for $q\in(0,d)$, and $q^*=dq/(d-q)$, $C^{S}_{d,q}$ is the best constant for the Sobolev's embeddings,
		\begin{equation*}
			\|\rho\|_{L^{q^*}}\leq C^{S}_{d,q}\|\nabla \rho\|_{L^q}.
		\end{equation*}
		The explicit value for these constants are known (see \cite{aubin_problemes_1976,talenti_best_1976,lieb_sharp_1983}).
	\end{remark}

	Note that for $d=2$, the point $(\ak,m)=(2,1)$ in the context of power law diffusion, $(\ak,\al)=(2,2)$ in the notations of fractional diffusion and $(\ak,p)=(2,2)$ in the notations of the present paper all correspond to the classical Keller-Segel equation, and the three different definitions of the fair competition case coincide for this equation.

\section{Proof of Theorem \ref{thm:main}}

	We begin this section by introducing the $p$-Fisher information $I_p$ on $\lt(W ^{1,p}\rt)^{p'} := \{\rho, \rho^{1/p'}\in W^{1,p}\}$ as
	\begin{equation*}
		I_p(\rho)=\intd \frac{\lt|\nabla \rho\rt|^p}{\rho}= (p')^p\|\nabla (\rho^{1/p'})\|^p_{L^p},
	\end{equation*}
	which is a generalization of the classical Fisher information (i.e. the case $p=2$). First remark that a straightforward computation using Hölder's and Sobolev's inequalities shows that $(W^{1,p})^{p'} \subsetArrow W^{1,p-1}_\loc$ so that $\Delta_p \rho$ is well-defined for $\rho$ with finite $p$-Fisher information. Then for any $p \in \lt(\frac{2d}{d+1},\frac{3d}{d+1}\rt)$, $q\in \lt[1,r\rt]$ with $r=\frac{p^*}{p'}$ and $\rho\in \lt(W ^{1,p}\rt)^{p'}\cap L^1$ it holds
	\begin{equation}\label{eq:GNS}
		\|\rho\|_{L^q} \leq ((p')^{-1}\CS_{d,p})^\frac{r'p'}{q'} \|\rho\|_{L^1}^{1-\frac{r'}{q'}} I_p(\rho)^\frac{r' p'}{q'p}.
	\end{equation}
	Indeed by Sobolev's embeddings it holds
	\begin{align*}
		\|\rho\|^{p/p'}_{L^r} = \|\rho^{1/p'}\|^p_{L^{p^*}} &\leq (C^S_{d,p})^{p}\|\nabla (\rho^{1/p'})\|^p_{L^p} \leq ((p')^{-1}C^S_{d,p})^{p}I_p(\rho),
	\end{align*}
	and using interpolation inequality between $L^1$ and $L^r$ yields the result.
	
	Then we need some tools in order to provide some moments estimate. First in the case $p\geq 2$ and  $k\in\lt[0,1\rt]$ there is $C> 0$ such that for any $\rho \in L^1 \cap \lt(W ^{1,p}\rt)^{p'}$
	\begin{align}\label{eq:momp>2}
		\intd (\Delta_p \rho) m &\leq C \|\rho\|^{\frac{p-1}{\al_p}}_{L^1} I_p(\rho)^\frac{\al_p-1}{\al_p}.
	\end{align}
	Indeed since $k\leq 1,$ by Hölder's inequality it holds
	\begin{align*}
		\intd \lt|\nabla \rho\rt|^{p-2}\nabla \rho \cdot \nabla m &\leq k\intd \rho^\frac{1}{p'} \rho^\frac{-1}{p'}\lt|\nabla \rho\rt|^{p-1} \weight{x}^{k-1}
		\\
		&\leq k\lt(\intd \rho^\frac{p}{p'} \rt)^\frac{1}{p} \lt(\intd \frac{|\nabla \rho|^p}{\rho}\rt)^\frac{1}{p'}=k\lt(\|\rho\|_{L^{\frac{p}{p'}}}I_p(\rho)\rt)^\frac{1}{p'}.
	\end{align*}
	Then, using inequality~\eqref{eq:GNS}, we obtain
	\begin{equation*}
		\|\rho\|_{L^{\frac{p}{p'}}}\leq C \|\rho\|_{L^1}^{\frac{1}{p'}\lt(1-\frac{r'(p-2)}{(p-1)}\rt)} I_p(\rho)^{\frac{r'p'}{\lt(p-1\rt)'p }},
	\end{equation*}
	and the result follows since 
	\begin{align*}
		\lt(\frac{r'p'}{\lt(p-1\rt)'p} + 1\rt)\frac{1}{p'} = \frac{\al_p-1}{\al_p}
		\text{ and }
		\frac{1}{p'}\lt(1-\frac{r'(p-2)}{(p-1)}\rt) = \frac{p-1}{\al_p}.
	\end{align*}
	Then in the case $p\in\lt(\frac{2d}{d+1},2\rt)$ and $k \in (0,\alpha_{p})$, there is $C> 0$ such that
	\begin{equation}\label{eq:momp<2}
		\intd (\Delta_p\rho) m \leq C \lt(\intd \rho m\rt)^\frac{1}{p'} I_p(\rho)^\frac{1}{p'}.
	\end{equation}
	Indeed by Hölder's inequality, since $p\leq 2$,
	\begin{align*}
		\intd \lt|\nabla \rho\rt|^{p-2}\nabla \rho \cdot \nabla m &\leq k\intd \lt|\nabla \rho\rt|^{p-1} \weight{x}^{k-1}
		\\
		&\leq k\intd \rho^\frac{1}{p'}\weight{x}^\frac{k}{p'} \rho^\frac{-1}{p'}\lt|\nabla \rho\rt|^{p-1} \weight{x}^{\frac{k}{p}-1}
		\\
		&\leq k\lt(\intd \rho m\rt)^\frac{1}{p'} \lt(\intd \frac{|\nabla \rho|^p}{\rho}\rt)^\frac{1}{p'} \lt(\intd \weight{x}^\frac{k-p}{2-p}\rt)^{\frac{2}{p} - 1},
	\end{align*}
	and the result follows since by assumption $\frac{k-p}{2-p} < -d$.

	\begin{demo}[ of Theorem~\ref{thm:main}]
		We only provide the a priori estimate necessary to the rigorous proof. Following the claim of \cite[Proof of Theorem 5.2, Step 1]{liu_degenerate_2016}, we can retrieve well posedness for the regularized problem \eqref{eq:agdif} with $K_\ak$ replaced with $K^\eps_\ak(x)=\indic_{|x|\geq \eps}K(x)+\indic_{|x|\leq \eps}\eps^{-\ak}x$. The preservation of positivity is a consequence of Kato's inequality for the $p$-Laplacian (see \cite{horiuchi_remarks_2001,liu_remarks_2016}). Then letting $\eps$ go to $0$ and using the a priori estimate we are about to prove with together with a standard compactness argument (similarly as what is done in \cite[Section 2.5]{blanchet_two-dimensional_2006})  provides the rigorous proof.
		
		\step{1. Entropy dissipation} We first estimate the dissipation of entropy using together Hardy-Littlewood-Sobolev's inequality and \eqref{eq:GNS} with $q=\frac{2d}{2d-\ak}$ as
		\begin{align*}
			\dt\intd\rho \log \rho
			&= -\intd \lt( \lt|\nabla \rho\rt|^{p-2} \nabla \rho \rt)\cdot \nabla \log \rho + \lambda \intd \divg((K_\ak *\rho)\rho)(\log \rho +1)
			\\
			&=-I_p(\rho) + \lambda \intd (\divg(K_\ak) *\rho)\rho
			\\
			&= -I_p(\rho) + \lambda (d-\ak) \iintd \frac{\rho(x)\rho(y)}{|x-y|^\ak}\d x\d y\\
			&\leq -I_p(\rho)+\lambda (d-\ak)\CHLS_{d,\ak,q} \lt(p'^{-1}\CS_{d,p}\rt)^{2\frac{r'p'}{q'}} M_0^{2-2\frac{r'}{q'}} I_p(\rho)^{2\frac{r' p'}{q' p}}.
		\end{align*}
		And since $2\frac{r' p'}{q' p}=\frac{\ak}{\al_p}$, $2-2\frac{r'}{q'}=2\lt(1-(p-1)\frac{\ak}{2\al_p}\rt)$ and $2\frac{r'p'}{q'}=p\frac{\ak}{\al_p}$, defining 
		\begin{equation*}
			C_{d,p}^{p-3}:= (d-\al_p){\CHLS_{d,\al_p,\frac{2d}{2d-\al_p}}} \lt(p'^{-1}\CS_{d,p}\rt)^{p},
		\end{equation*}
		we conclude this step with
		\begin{align}\label{eq:disentFC}
			\dt \intd \rho \log \rho &\leq
			\begin{cases}
			-\lt(1- \lambda C_{d,p}^{p-3}M_0^{3-p}\rt) I_p(\rho), \text{ if } \ak=\al_p
			\\ 
			-\frac{1}{2} I_p(\rho)+C \text{ if }\ak<\al_p.
			\end{cases}
		\end{align}
				
		\step{2. Moment estimate} 
		First in the case $p\geq 2$, we choose $k\in ([1-\ak]_+,1)$ and use \eqref{eq:momp>2}, symmetry and Young's inequality for any $\eps>0$ to obtain
		\begin{align*}
			\dt\intd \rho m
			&\leq C M_0^\frac{p-1}{\al_p} I_p(\rho)^\frac{\al_p-1}{\al_p} - \frac{\lambda}{2} \iintd K_\ak(x-y)\cdot \lt(\nabla m(x)-\nabla m(y)\rt)  \rho( \dd x)\rho(\dd y)
			\\
			&\leq C_\eps M_0^{\frac{p-1}{\al_p} \lt(\frac{\al_p}{\al_p-1}\rt)' } + \eps I_p(\rho)- \frac{\lambda}{2} \iintd K_\ak(x-y)\cdot \lt(\nabla m(x)-\nabla m(y)\rt)  \rho(\dd x)\rho(\dd y).
		\end{align*}
		Then in the case $p\in\lt(\frac{2d}{d+1},2\rt)$, we choose $k\in ([1-\ak]_+,\al_p)$, use \eqref{eq:momp<2} and obtain
		\begin{align*}
			\dt\intd \rho m &\leq C \lt(\intd \rho m\rt)^\frac{1}{p'} I_p(\rho)^\frac{1}{p'} - \frac{\lambda}{2} \iintd K_\ak(x-y)\cdot \lt(\nabla m(x)-\nabla m(y)\rt)  \rho(\dd x)\rho(\dd y) 
			\\
			&\leq C_\eps \lt(\intd \rho m \rt)^\frac{p}{p'}+ \eps I_p(\rho)- \frac{\lambda}{2} \iintd K_\ak(x-y)\cdot \lt(\nabla m(x)-\nabla m(y)\rt)  \rho(\dd x)\rho(\dd y).
		\end{align*}
		The last term in the r.h.s is dealt similarly as in \cite[Proof of Proposition 3.1]{lafleche_fractional_2018} and in any case, we end up with 
		\begin{equation}\label{eq:mom}	
		\begin{aligned}
			\dt\intd \rho m \leq \eps I_p(\rho)+C\lt(1 + \intd \rho m\rt),
		\end{aligned}
		\end{equation}
		where $C$ only depends on $d,p,\ak,\lambda,\eps,k$ and $M_0$
		
		\step{3.Conclusion} 
		We will now only treat the case $\ak=\al_p$ and 
		$\lambda M_0^{3-p}C_{d,p}^{p-3}< 1$, since the case $\ak<\al_p$ can be treated even more straightforwardly. For $k\geq 0$ denote $\nu_k > 0$ such that $\intd e^{-\nu_k m(x)}\d x =1$, and recall that, with $h(u)=u\ln u -u+1\geq 0$, it holds
		\begin{align*}
			\intd \frac{\rho}{M_0} \ln \frac{\rho}{M_0} &= \intd h\lt(\frac{\rho}{M_0} e^{\nu_k m}\rt)e^{-\nu_km} + \intd \frac{\rho}{M_0} \ln (e^{-\nu_km})\geq -\nu_k \intd \frac{\rho}{M_0}m,
		\end{align*}
		and then 
		\begin{equation*}
			\intd \rho\ln \rho \geq M_0\ln M_0-\nu_k \intd \rho m,
		\end{equation*}
		which yields for fixed $\nu > \nu_k$, combining linearly \eqref{eq:disentFC} and \eqref{eq:mom}
		\begin{align*}
			\lt(\nu - \nu_k\rt)\intd \rho m &\leq - M_0\ln M_0
			+ \intd\rho \log \rho + \nu\intd \rho m
			\\
			&\leq \intd\rho^\init \log \rho^\init + \nu\intd \rho^\init m
			+ \nu C \int_0^t\lt(\intd \rho(s) m +1\rt) \d s
			\\
			&\quad - \lt(1- \lambda C_{d,p}^{p-3} M_0^{3-p} - \eps \nu\rt)\int_0^t I_p(\rho)(s) \d s.
		\end{align*}
		Therefore, for $\eps > 0$ small enough, $\rho\in L^\infty_\loc\lt(\R_+,L^1_k\rt)$ by Gronwall's inequality. We emphasize that this estimate also applies to the $p$-heat equation, i.e. \eqref{eq:pLapKS} with $\lambda=0$. Finally coming back to \eqref{eq:disentFC} yields
		\begin{align*}
			\lt(1- \lambda C_{d,p}^{p-3} M_0^{3-p} \rt)\int_0^t I_p(\rho)(s)\d s &\leq \intd \rho^\init \log \rho^\init - \intd \rho \log \rho
			\\
			&\leq \intd \rho^\init \log \rho^\init + \nu_k \intd \rho m-M_0\ln M_0,
		\end{align*}
		and we conclude the proof using inequality~\eqref{eq:GNS}.
	\end{demo}
	
\section*{Acknowledgements}	

The second author was supported by the Fondation des Sciences Math\'ematiques de Paris and Paris Sciences \& Lettres Universit\'e

%% ********************  Bibliographie  ********************

\bibliographystyle{abbrv} % apalike, ieee, plain, alpha, unsrt, abbrv
\bibliography{KS_p-lap}

\bigskip
\signll
\signss

\end{document}